\documentclass[final]{elsart}
\usepackage{epsfig}
\usepackage{amsmath}
\usepackage{amssymb}

\newcommand{\imm}{\mathrm{i}}
\newcommand{\diff}{\mathrm{d}}

\newcommand{\N}{\mathbb{N}}
\newcommand{\Z}{\mathbb{Z}}
\newcommand{\R}{\mathbb{R}}
\newcommand{\C}{\mathbb{C}}

\newcommand{\toe}{T}
\newcommand{\cir}{\mathcal{C}}

\newcommand{\D}{\displaystyle}
\newcommand{\mirror}{\mathcal{M}}
\newcommand{\diag}{{\rm diag}}
\newcommand{\norm}[2][]{{\|#2\|}_{#1}}

\renewcommand{\v}[1]{\ensuremath{\mathbf{#1}}}

\newtheorem{proposition}{Proposition}
\newtheorem{remark}{Remark}
\newtheorem{definition}{Definition}

\newenvironment{proof}
{\emph{Proof.}} {\hspace*{\fill}$\Box$}

\newcounter{ChToep}
\newenvironment{prolist}
{\begin{list}{(\textbf{\Alph{ChToep}})}{\usecounter{ChToep}}}
{\end{list}}

\begin{document}
\begin{frontmatter}

\title{A note on grid transfer operators for multigrid methods\thanksref{titlefn}}
\author{Marco Donatelli}
\address{Dipartimento di Fisica e Matematica, Universit{\`a} dell'Insubria - Sede di Como,
    Via Valleggio 11, 22100 Como, Italy.}
\thanks[titlefn]{Supported by MUR grant number 2006017542.}
\ead{marco.donatelli@uninsubria.it}
\ead[url]{http://scienze-como.uninsubria.it/mdonatelli}

\maketitle
%
\begin{abstract}
The Local Fourier analysis (LFA) is a classic tool to prove
convergence theorems for multigrid methods (MGMs). In particular, we
are interested in optimality that is a convergence speed independent
of the size of the involved matrices. For elliptic partial
differential equations (PDEs), a well known optimality result
requires that the sum of the orders of the grid transfer operators
is not lower than the order of the PDE to solve. Analogously, when
dealing with MGMs for Toeplitz matrices in the literature an
optimality condition on the position and on the order of the zeros
of the symbols of the grid transfer operators has been found. In
this work we show that in the case of elliptic PDEs with constant
coefficients, the two different approaches lead to an equivalent
condition. We argue that the analysis for Toeplitz matrices is an
algebraic generalization of the LFA, which allows to deal not only
with differential problems but also for instance with integral
problems. The equivalence of the two approaches gives the
possibility of using grid transfer operators with different orders
also for MGMs for Toeplitz matrices. We give also a class of grid
transfer operators related to the B-spline's refinement equation and
we study their geometric properties. This analysis suggests further
links between wavelets and multigrid methods. A numerical
experimentation confirms the correctness of the proposed analysis.
\end{abstract}
\end{frontmatter}
%
%
\section{Introduction}\label{sect:intr}
Multigrid methods (MGMs) are widely used for solving elliptic PDEs.
The convergence analysis is usually done in the case of constant
coefficients. Let us consider standard finite differences
discretization for the following $d$-dimensional problem
\begin{equation} \label{eq:pde}
    \left\{
    \begin{array}{l}\displaystyle
        (-1)^q \sum_{i=1}^d \frac{\diff^{2q}}{\diff x_i^{2q}}u(x)
        = g(x), \qquad x \in \Omega=(0,1)^d, \; q \geq 1,\\
        \mbox{boundary conditions on} \; \partial\Omega,
    \end{array}\right.
\end{equation}
where $x=(x_1, \dots,x_d)$.
For simplicity of the presentation we assume the same meshsize $h$
for each dimension, but generalization to $h=(h_1, \ldots, h_d)$ is possible.
Hence, approximating \eqref{eq:pde} by
centered finite differences on a uniform grid of $n$ subintervals of size $h$
in each dimension, we obtain the linear system
\begin{equation}\label{eq:sisdiscr}
    A_n\v{y}=\v{b}
\end{equation}
of order $n^d \times n^d$.
Neglecting the boundary conditions, the matrix $A_n$ is
a $d$-level Toeplitz matrix and it is banded at each level.

Let $L_h$ be a discretization of the differential
operator in \eqref{eq:pde}, then its Fourier
transform is
\begin{equation}\label{eq:Lhat}
    \hat{L}(\omega) = \sum_{j\in\Z^d}{l_j\e^{\imm \langle jh | \omega
    \rangle}},
\end{equation}
for $\omega \in [-\pi/h, \pi/h]^d$, $\imm^2=-1$, and
\begin{equation}\label{eq:Lcoeff}
    l_j = \frac{h^d}{(2\pi)^d} \int_{[-\pi/h, \pi/h]^d} \hat{L}(\omega) \e^{-\imm  \langle jh | \omega \rangle} \diff\omega,
\end{equation}
where the operations between multi-indices are intended component
wise and $\langle\,\cdot\,|\,\cdot\,\rangle$ denotes the usual
scalar product between vectors. Using centered finite differences of
precision $2$ and minimal bandwidth, the polynomial
$\hat{L}(\omega)$ has degree $2q$, that is the order of the PDE
\eqref{eq:pde}, and $L_h$ is completely defined from its
$d$-dimensional stencil formed by the coefficients $l_j$, with $j_i
= -2q+1,\,  \dots, \, 2q-1$ for $i=1,\dots,d$.

It is well-known that MGMs are optimal solvers
for PDEs of the form \eqref{eq:pde}, i.e., they require about a
constant number of iterations varying $h$ and each iteration has an
arithmetic cost proportional to the matrix-vector product
\cite{Trot}. Such property is obtained imposing a well known
condition on the order of the grid transfer operators:
\begin{equation}\label{cond:ordmgm}
m_r + m_p \geq m,
\end{equation}
where $m=2q$ is the order of the PDE, $m_r$ is the order of the
restriction and $m_p$ is the order of the prolongation
\cite{LFA_Brandt}, denoted as high frequencies order in
\cite{Hemker,Yav}. The condition (\ref{cond:ordmgm}) follows from
the LFA for the two grid method (TGM). In order to obtain more
powerful grid transfer operators, that is in order to devise an
optimal MGM, inequality (\ref{cond:ordmgm}) should to be satisfied
strictly \cite{Trot}. We note that the LFA does not consider the
border effects, i.e., it assumes periodic boundary conditions or an
infinite domain \cite{LFA_Brandt}.

MGMs for multilevel positive definite Toeplitz matrices have been
developed in the years looking only to the linear system
\eqref{eq:sisdiscr}, independently of the continuous problem
\cite{FS1,FS2,SConv,ADS,Chan,Huk}. A MGM for Toeplitz matrices was
early defined in \cite{FS1} using a powerful eigenvalue interlacing
property with the matrices in the $\tau$ algebra (the class of
matrices diagonalized by discrete sine transforms of type I). This
first proposal was extended to the multilevel case in
\cite{FS2,Chan,SConv}. Since Toeplitz matrices do not define an
algebra and hence are difficult to manipulate, convergence results
are provided using matrix algebra approximations like $\tau$ or
circulant matrices having the same spectral distribution of the
Toeplitz matrices. In other words, we require that the circulant or
the $\tau$ approximates share the same symbol of the original
Toeplitz matrix \cite{DiBen,T2}. In this paper we consider the
circulant case. A MGM for circulant matrices was introduced in
\cite{Serra-Possio}. Furthermore in \cite{ADS, ADmultid},
generalizing the techniques used in \cite{SConv}, and using the Ruge
and St{\"u}ben theory \cite{RStub} and the Perron-Frobenius theorem,
a complete proof of the optimality of the $V$-cycle for multilevel
circulant and $\tau$ matrices was proposed. This analysis leads to a
stronger condition with respect to the previous two grid analysis.

In this paper we show that the techniques used in
\cite{FS1,FS2,SConv,Serra-Possio} represent a linear algebra
generalization of the LFA \cite{LFA_Brandt,Hemker}, in the case of
the Galerkin approach. Indeed, they lead to a condition analogous to
(\ref{cond:ordmgm}), but on the order of the zeros of the generating
functions of the grid transfer operators. The letter represents a
wide generalization since the case of discretization of elliptic
PDEs of order $2q$ corresponds to the case of generating function
which are nonnegative (ellipticity), with a unique zero at zero of
order $2q$ (consistency condition). In other words, the case of
discretization of elliptic PDEs is a subcase of nonnegative symbols
with unique zero at zero, which in turn represents the case when the
algebraic problem is ill-conditioned in a subspace of low
frequencies. Therefore by using the Toeplitz approach other cases
can be considered including the case when the ill-conditioning
arises in high-frequencies: we recall that the latter characterizes
some integral problem related to signal/images restoration. We will
show that considering the problem \eqref{eq:pde}, the LFA done in
\cite{LFA_Brandt,Hemker,Yav} and the Toeplitz approach
\cite{FS1,SConv} (which was introduced independently) are
essentially equivalent. As already stressed, second approach is more
general since it can be applied also when the zero of the symbol is
not at the origin, or there exist several zeros (multiple sources of
ill-conditioning). By the Galerkin approach, we have the only
limitation that the restriction must be proportional to the
transpose of the prolongation, but in this paper we will show that
in practical implementations this condition is not necessary. More
precisely, we will generalize the MGM for Toeplitz matrices to the
case of a restriction different to the transpose of the
prolongation. A first suggestion to consider the linear algebra
tools for Toeplitz matrices as a generalization of the LFA for
multigrid methods was given in \cite{HukFour,GLT}.

In this paper we also define a class of grid transfer operators that
satisfy the conditions in \cite{Hemker} but that are not
interpolating.
More in detail the
considered operators are defined looking for the smallest support of
the symbol for a fixed order and are related to the refinement
equation of B-spline. We will give a geometrical interpretation of
the operator of order 4 and it will be compared with the cubic
interpolation. These B-spline grid transfer operators allow us to
discuss some relations between wavelets and multigrid methods.
Eventually, a numerical experimentation validates our proposals.

The paper is organized as follows. Firstly in {\S}\ref{sect:LFA} we
present the results in \cite{Hemker}. In {\S}\ref{sect:MGM_toep} we
describe the MGM defined in \cite{FS1, FS2, SConv} for multilevel
Toeplitz matrices using the zeros of the generating functions and we
compare the two-grid analysis with the results in \cite{Hemker}. In
{\S}\ref{sec:gentoep} we generalize the TGM described in
{\S}\ref{sect:MGM_toep} to the case of restriction not necessary
proportional to the transpose of the prolongation. In
{\S}\ref{sec:b-spline} we give a new class of grid transfer
operators with minimum support for a fixed order and we show that
such class is related to the B-spline. This allows, in
{\S}\ref{sec:wave}, some observations about the relations between
wavelet and multigrid methods. In {\S}\ref{sect:exp} some numerical
results validate the previous proposals both for Toeplitz
non-differential problems and for PDEs with nonconstant
coefficients. The final, {\S}\ref{sect:concl} is devoted to some
concluding remarks.

\section{The low and high frequencies order analysis}\label{sect:LFA}
We introduce a grid transfer operator that is not effective alone
(it has order zero), but which is the basic tool for developing more
powerful projectors. It is the classic \emph{down-sampling}
operator, called \emph{elementary restriction} in \cite{Hemker} and
\emph{cutting matrix} in \cite{FS1}. In the one-dimensional case, we
set $n = n^{(0)} > n^{(1)} > \dots > n^{(l)} > 0$, $l \in \N$, such
that $n^{(i+1)} = (n^{(i)}-(n^{(i)} \, \mathrm{mod} \, 2))/2$ and we
define the down-sampling matrix $K_{n^{(i)}} \in \R^{n^{(i+1)}
\times n^{(i)}}$ as
\begin{equation}\label{eq:cut}
    [K_{n^{(i)}}]_{j,k} =
    \left\{
    \begin{array}{l l}
        1 & \mbox{if } j = 2k - (n^{(i)}+1) \, \mathrm{mod} \, 2, \\
        0 & \mbox{otherwise,}
    \end{array}
    \right.
    \qquad
    k = 1,\dots,n^{(i+1)}.
\end{equation}
In the $d$-dimensional case the down-sampling matrix is defined by tensor product as
$K_{n^{(i)}} = K_{n^{(i)}_1} \otimes K_{n^{(i)}_2} \otimes \dots \otimes K_{n^{(i)}_d}$.

Higher order grid transfer operators are defined by convolution
with the down-sampling operator. The prolongation is
$P_{n^{(i)}}(p_i)=T_{n^{(i)}}(p_i)K_{n^{(i)}}^T$, while the
restriction is $R_{n^{(i)}}(r_i)=K_{n^{(i)}}T_{n^{(i)}}(r_i)$
that usually is the transpose, up to a constant factor, of the
prolongation.
The matrix $T_n(t)$ is the $d$-level Toeplitz matrix generated by the function $t$.
The Toeplitz matrices will be described in Section \ref{sect:MGM_toep}.
The symbols $r_i$ and $p_i$ should be trigonometric
polynomials of low order to maintain the computational cost of the
matrix vector product proportional to $O(N(n^{(i)}))$.

Since more powerful grid transfer operators give a greater
computational cost, it is important to find sufficient conditions
such that we can decide for a fixed problem the cheapest grid
transfer operators that allows to obtain an optimal MGM. This task
can be done using the LFA \cite{LFA_Brandt,Hemker} obtaining the
condition \eqref{cond:ordmgm}.

We define the set of all corners of $x$ as $\Omega(x) = \{\, y \, |
\,  y_j \in \{x_j, \, \pi + x_j\}, \, j=1, \dots, d \}$. With the
change of variable $x=\omega h$, the set of all frequencies on the
fine grid that correspond to the frequency $\omega$ on the coarse
grid is $\{z=y/h \; | \; y\in\Omega(x) \}$. Moreover, according to
the terminology in \cite{FS1} we define the set of the ``mirror''
points of $x$ as $\mirror(x)=\Omega(x)\setminus \{x\}$. Since in the
rest of the section we will consider only two grids, $n$ will denote
the fine grid. Moreover, for unifying the treatment, a generic grid
transfer operator is denoted by $B_n(g)$ where $g$ is multiplied by
a factor $2^d$ when $B_n(g)$ is the prolongation, i.e., $B_n(g) =
R_n(g)$ or $B_n(g) = P_n(2^dg)$.
\begin{definition}\label{def:lf}
    The \emph{Low Frequency order (LF)} of a grid transfer operator $B_n(g)$
    is the largest number $s \geq 0$ for which
    \[g(x) = 1 + O(|x|^s), \qquad \mbox{ for } |x| \to 0.\]
\end{definition}
\begin{definition}\label{def:hf}
    The \emph{High Frequency order (HF)} of a grid transfer operator $B_n(g)$ is
    the largest number $s \geq 0$ for which
    \[g(y) = O(|x|^s), \qquad \forall y \in \mirror(x), \; \mbox{ for } |x| \to 0.\]
\end{definition}
For $x=\omega h$, $|x|\rightarrow 0$ means $h \rightarrow 0$ since $\omega$ is fixed.
For $h=(h_1,\dots,h_d)$ we can define $|x| = \max_{i=1,\dots,d}(|x_i|)$.

For the grid transfer operators LF and HF are more
general then classic interpolation order.
\begin{proposition}[\cite{Hemker}]\label{prop:ordint}
  \begin{itemize}
    \item[(i)] If a restriction leaves all polynomials of degree $s-1$ invariant,
            then the LF of the operator is $s$.
    \item[(ii)] If a prolongation leaves all polynomial of degree $s-1$ invariant,
            then both the LF and HF are at least $s$.
  \end{itemize}
\end{proposition}
For instance the linear interpolation has LF = HF = 2,
while the cubic interpolation has LF = HF = 4.

Furthermore, we can derive the condition (\ref{cond:ordmgm}) from the following
\begin{proposition}[\cite{Hemker}]\label{prop:hemk}
    Given a constant-coefficient, linear differential operator of order $m$,
    a necessary condition for non-increasing the high frequencies arising
    from a coarse grid correction with two grids it is
    \begin{equation}\label{eq:condF}
        \gamma_r + \gamma_p \geq m,
    \end{equation}
    where $\gamma_p$ and $\gamma_r$ are the HF of the prolongation and of the
    restriction respectively.
\end{proposition}

From Proposition \ref{prop:ordint} part \textit{(ii)}  the condition
\eqref{eq:condF} is a generalization of the analogous condition on
the interpolation order. The LF is important for the restriction
thanks to Proposition \ref{prop:ordint} part \textit{(i)}, but it
seems not necessary for the two grid analysis in Proposition
\ref{prop:hemk}. However, in \cite{LFA_Brandt} it is shown that for
an efficient MGM a further condition is that both LF and HF are
positive. This further request arises also from the Galerkin
approach (see \cite{Yav}) and is natural for obtaining an effective
MGM.

Eventually, we note that, since the grid transfer operation has to
be computationally cheap, the function $g$ in Definitions
\ref{def:lf} and \ref{def:hf} should be a trigonometric polynomial
of low degree. Moreover, from Proposition \ref{prop:ordint} a good
class of grid transfer operators should have at least LF$>0$.
Interpolating operators define a class with LF=HF. A further class
of operators with a fixed HF and LF$>0$ will be described in
Section~\ref{sec:b-spline}.

\section{A MGM for Toeplitz matrices by generating functions}\label{sect:MGM_toep}
In this section we briefly introduce the MGM defined and analyzed in
\cite{FS1, FS2, Chan, TGM-wave, SConv, Serra-Possio, ADS, ADmultid}
for the multidimensional $\tau$, circulant, Toeplitz and other
matrix-algebras related to trigonometric transforms.

Toeplitz matrices arise from the discretization of convolution
operators with a shift invariant kernel and hence not only from
PDEs, but also from several other applications, e.g., image
deblurring problems \cite{ADS}. Toeplitz matrices are completely
defined by the matrix size and the symbol also called generating
function. Let $f$ be a continuous function on $\R^d$ and having
period $2\pi$ with respect to each variable, the Fourier
coefficients of $f$ are defined as
\begin{equation}\label{a_j=}
  a_j = \frac1{(2\pi)^d}
        \int_{{[-\pi,\pi]}^d} f(x)e^{-\imm \langle j|x\rangle}\ dx,
  \qquad  j\in\Z^d.
\end{equation}

\begin{remark}\label{rem0}
With the change of variable $x=\omega h$, it holds $a_j=l_j$ and
$f(x) =\hat{L}(\omega)$.
\end{remark}

From the coefficients $\{a_j\}$ one can build \cite{T2} the sequence
$\{\toe_n(f)\}$ of multilevel Toeplitz matrices. Every matrix
$\toe_n(f)$ is explicitly written as
\[
  \toe_n(f) = \sum_{|j_1|\leq  n-1}\ \dots\ \sum_{|j_d|\leq n-1}
              a_{(j_1,\dots,j_d)}
              J_{n}^{[j_1]}\otimes\dots\otimes J_{n}^{[j_d]}.
\]
Here $\otimes$ denotes the usual tensor product and $J^{[j_i]}_n \in
\R^{n \times n}$ is the matrix whose entry ($s,t$) equals $1$ if
$s-t=j_i$ and is $0$ elsewhere, for $i=1,\dots,d$. Many structural
and spectral properties of $\toe_n(f)$ derive from its generating
function $f$. Indeed, if $f$ is real valued, then $a_{-j}=\bar{a}_j$
for every $j$ and the matrices $\toe_n(f)$ are Hermitian for every
$n$; if $f$ is also non-negative but not identically zero then
$\toe_n(f)$ is positive definite.

\begin{remark}
The main difference between $f(x)$  and $\hat{L}(\omega)$ is that
$\omega$ denotes the frequency for the current discretization step
$h$, information that seems to be lost in $f$, but that comes out
from the matrix $T_n(f)$ regarding the current discretization
($h=1/(n+1)$). For instance, let $L_h$ be the three point
discretization of the Laplacian: then $l_0 = 2/h^2$ and $l_{-1} =
l_{1} = -1/h^2$. On the other hand, in the algebraic approach for
Toeplitz matrices the constant factor $1/h^2$ is moved to the right
hand side (rhs) obtaining $a_0=2$ and $a_{-1} = a_1 = -1$. However
the information of the order 2 of the Laplace operator is preserved
since $f(x)=2-2\cos(x)$ vanishes at the origin with order 2. More in
general, discretizing \eqref{eq:pde} with finite centered
differences of minimal precision and moving the coefficient
$1/h^{2q}$ to the rhs, by consistency, the symbol $f(x)$ vanishes at
the origin with order $2q$.
\end{remark}

Convergence results for MGMs for PDEs are usually obtained
neglecting the boundary conditions. In a similar way, MGMs for
Toeplitz matrices are defined starting from matrix algebras like
$\tau$ or circulant. Imposing periodic boundary conditions in
\eqref{eq:pde}, the matrix $A_n$ in \eqref{eq:sisdiscr} is
circulant. Circulant matrices are simultaneously diagonalized by the
Fourier transform $F_n = \frac{1}{\sqrt{n}}[e^{-\imm j
y^{(n)}_i}]_{i,j}$, where $y^{(n)}_i=2\pi i/n$, $i=0,\dots,n-1$.
More precisely, the algebra of the circulant matrices can be
formally defined as $  \left\{A_n \, | \, A_n=F_n\cdot\diag{(\v{d})}
\cdot F^{H}_n, \; \v{d}\in\C^n \right\} $ where the vector $\v{d}$
of the eigenvalues is equal to $f(\v{y}^{(n)})$, $\v{y}^{(n)} =
(y^{(n)}_0, \dots, y^{(n)}_{n-1})$, and a circulant matrix will be
denoted by $\cir_n(f)$. In the $d$-dimensional case the indices
involved are multiindices, $F_n = F_{n_1} \otimes\dots\otimes
F_{n_d}$ has size $N(n)$ and
$\v{y}^{(n)}=\v{y}^{(n_1)}\times\dots\times\v{y}^{(n_d)}$, where
$\times$ denotes the cartesian product.


We do not consider boundary effects, thus we will discuss only the
circulant case assuming periodic boundary conditions. In such case,
in order to maintain the same circulant structure at each level, we
have to start with $n=n^{(0)}=2^a$, where $a \in \N^d$. Moreover the
grid transfer operators are defined as
$P_{n^{(i)}}(p_i)=C_{n^{(i)}}(p_i)K_{n^{(i)}}^T$ and
$R_{n^{(i)}}(r_i)=K_{n^{(i)}}C_{n^{(i)}}(r_i)$. In our case
$A_n=\cir_n(f)$ is singular since $f$ vanishes at the origin which
is a grid point. However, without losing generality, we assume $A_n$
nonsingular replacing $f$ with its stabilized version that is by
correcting $A_n$ by adding a special rank-one matrix. This
correction is not consider here since it does not imply particular
assumptions but it leads only to unnecessary complications in the
notation \nolinebreak \cite{ADS}.

Using the Galerkin approach, we must have $R_{n}(r)=P_{n}(p)^H$,
\i.e., $r=p$, and $A_{n/2} = P_{n}(p)^H A_{n} P_{n}(p)$. Thanks to
the structure of $P_{n}(p)$ we obtain that $A_{n/2}$ belongs again
to the circulant algebra \cite{Serra-Possio}. Thanks to the
Ruge-St{\"u}ben theory \cite{RStub}, the TGM and the V-cycle
convergence analysis can be split in two independent conditions, one
on the smoother and the other on $p_i$, for $i=0,\dots,l-1$, i.e.,
on the grid transfer operators.

\begin{remark}
Several simple iterative methods, like relaxed Jacobi,
satisfy the smoothing condition, therefore the main task is the
study of the approximation condition for the grid transfer operators.
\end{remark}

In \cite{Serra-Possio} the optimality of the TGM was proved for
circulant matrices, under the following conditions on the grid
transfer operators.

\begin{proposition}[\cite{Serra-Possio}]\label{TGM_Serra}
Let the coefficient matrix be $A_n = \cir_n(f)$ with $f$ having a
unique zero at $x^0$. Defining $P_n(p) = \cir_n(p)K_n^T$ and
$R_n(r)=\alpha P_n(p)^H$, i.e. $r=\alpha p$, $\alpha \in \R
\backslash \{0\}$, where $p$ is a trigonometric polynomial non
identically zero and such that for each $x \in [-\pi, \pi)^d$
\begin{subequations}\label{ipo}
    \begin{equation}\label{ipo_a}
        \limsup_{x\to x^0}\left|\frac{p(y)^2}{f(x)}\right| = c < +\infty,
        \qquad \forall \, y \in \mirror(x),
    \end{equation}
    \emph{where}
    \begin{equation}\label{ipo_b}
        \sum_{y\in \Omega(x)} p(y)^2 > 0,
    \end{equation}
\end{subequations}
then defining $A_{n/2} =  \alpha P_{n}(p)^H A_{n} P_{n}(p)$ the TGM is optimal.
\end{proposition}

\begin{proof}
For $\alpha=1$ see \cite{Serra-Possio}.
For $\alpha \neq 1$ it is enough to observe that the coarse grid
correction $CGC = I-P_n(R_nA_nP_n)^{-1}R_nA_n$ is independent of $\alpha$.
\end{proof}

In order to compare this result with the Proposition \ref{prop:hemk}
we have to require $p=2^d r$, thus the \eqref{eq:condF} becomes
$2\gamma_r \geq m$. We show the equivalence between two different
convergence analysis for elliptic PDEs with constant coefficients:
the LFA described in Section \ref{sect:LFA} and the analysis for
Toeplitz matrices based on the zeros of the generating functions
described here.

\begin{proposition}\label{prop:equiv}
Let $P_n(p)=R_n(2^d r)^H$, discretizing \eqref{eq:pde} by finite
centered differences of order 2 and minimal bandwidth, the
conditions \eqref{eq:condF} and \eqref{ipo_a} are equivalent.
\end{proposition}
\begin{proof}
By Definition \ref{def:hf} $B_n(g)$ has HF$=s$ if and only if
$g(y)=0$ with order $s$ for all $y \in \mirror(x)$. The
discretization of an elliptic constant coefficient PDE of order $m$
by finite centered differences of precision 2 and minimal bandwidth
leads to $A_n=\cir_n(f)$ (in the case of periodic BCs) with $f$
vanishing at the origin with order $m$ (see Remark \ref{rem0}). From
condition (\ref{ipo_a}) $p$ (or equivalently $r$) must be chosen
such that $p(y)=0$, for all $y \in \mathcal{M}(0)$ with order
$2\gamma_p \geq m$. This is exactly the inequality in
\eqref{eq:condF}.
\end{proof}

The previous proposition shows that in the case of $p=2^dr$ and
using the Galerkin approach, condition \eqref{ipo_a} is a
generalization of condition \eqref{eq:condF} to general problems not
necessarily of differential type. The main difference between the
two approaches relies in the coarse strategy. The results in
\cite{Hemker} and summarized in Section \ref{sect:LFA} assume a
discretization of the same PDE with the same formula at each grid.
This imposes a right scaling of the grid transfer operators (i.e.,
$B_n(g)$ has LF $>0$ iff $g(0)=1$). The latter is not necessary in
the Galerkin approach adopted by the Toeplitz analysis, since the
coarse matrix is defined as $A_{n/2} =  \alpha P_{n}(p)^H A_{n}
P_{n}(p)$. Indeed the condition \eqref{ipo_b} requires only $p(0)
\neq 0$. More specifically, $p$ can be defined up to a scaling
factor since this gives only a different scaling of $A_{n/2}$.
However, the two approaches are comparable because from item 2 in
Proposition \ref{lowlev1d}, the coarse problem vanishes again at the
origin and with the same order of the finer problem \cite{SConv}.
Using the PDE language, this means that for the Galerkin approach
the linear system at the coarse grid is essentially (neglecting
boundary conditions) the discretization of the same PDE with a
formula of the same order.

At the end of Section \ref{sect:LFA} we noted that Proposition
\ref{prop:hemk} does not requires any condition on the LF of the
grid transfer operators. The only interest on the LF could be
deduced from Proposition \ref{prop:ordint}, and mainly for the
restriction. On the other hand, the TGM condition \eqref{ipo_b}
requires that the grid transfer operators have a positive LF (up to
a scaling factor). This is exactly the same requirement obtained in
\cite{Yav} for the Galerkin strategy and in \cite{LFA_Brandt} for an
efficient MGM. In fact a condition LF$ =0$ is equivalent to violate
\eqref{ipo_b} which implies $p(x^0)=0$. As a consequence, the
associated grid transfer operators could fail to be full rank. The
latter produces an increase of the ill-conditioning and could lead
to singularity at the lower levels with a potential substantial
change in the subspace related to small eigenvalues.

\section{MGM for Toeplitz matrices with a prolongation different
from the transpose of the restriction}\label{sec:gentoep}
In practical implementation the
condition
\begin{equation}\label{eq:p=r}
    R_n(r)=\alpha P_n(p)^H
\end{equation}
seems to be not necessary. The only request is that $A_{n/2} =
R_{n}(r) A_{n} P_{n}(p)$ is again positive definite for a recursive
application of the algorithm. On the other hand, the condition
\eqref{eq:p=r} is very useful for a theoretical analysis, because if
$r \neq p$ the coarse grid correction $CGC =
I_n-P_n(p)A_{n/2}^{-1}R_n(r)A_n$ is again a projector, but it is not
longer unitary with respect to the scalar product
$<\!\v{y},\,\v{x}\!>_{A_n}=\v{y}^H A_n \v{x}$, $A_n$ Hermitian
positive definite, for all $\v{y}$, $\v{x} \in \C^{n}$. For the well
definiteness of a MGM, mainly to ensure that the same smoother is
convergent also to the coarse levels, $A_{n/2}$ should be positive
definite to apply recursively the algorithm. This condition is easy
satisfied for $p_i \geqslant 0$, $r_i \geqslant 0$ (not identically
zero) and either both even or both odd, $i=0,\dots,l-1$. More
generally we could use $p_ir_i\geq 0$ with isolated zeros,
$i=0,\dots,l-1$. Therefore, the following generalization of
Proposition \ref{TGM_Serra} can be conjectured.

\textsc{TGM conditions}. \emph{Let the coefficient matrix be $A_n =
\cir_n(f)$, with $f$ having a unique zero at $x^0$. Defining $R_n =
K_n\cir_n(r)$ and $P_n = \cir_n(p)K_n^T$ where $p$ and $r$ are
trigonometric polynomials non identically zero and such that for
each $x \in [-\pi, \pi)^d$}
\begin{subequations}\label{NEWipo}
    \begin{equation}\label{NEWipo_a}
        \limsup_{x\to x^0}\left|\frac{r(y)p(y)}{f(x)}\right| = c < +\infty,
        \qquad \forall \, y \in \mirror(x),
    \end{equation}
    \emph{where}
    \begin{equation}\label{NEWipo_b}
        \sum_{y\in \Omega(x)} r(y)p(y) \neq 0.
    \end{equation}
\end{subequations}
\emph{Then, defining $A_{n/2} =  R_{n}(r) A_{n} P_{n}(p)$, the TGM
is optimal}.

These two conditions are motivated by the analysis in the previous
section and by the following Proposition \ref{prop:equiv2} that
extends Proposition \ref{prop:equiv} to the case of $r\neq\alpha p$.
Moreover, the numerical experiments in Section \ref{sect:exp} will
validate these conditions.

\begin{proposition}\label{prop:equiv2}
Discretizing \eqref{eq:pde} by finite centered differences of order
2 and minimal bandwidth, the conditions \eqref{eq:condF} and
\eqref{NEWipo_a} are equivalent.
\end{proposition}
\begin{proof}
The proof is analogous to that of Proposition \ref{prop:equiv}. It
is enough to observe that if $r$ and $p$ vanish at $y$ with order
$\gamma_r$ and $\gamma_p$ respectively, then $rp$ vanishes at $y$
with order $\gamma_r+\gamma_p$.
\end{proof}

We provide a further result useful to implement the corresponding MGM.

\begin{proposition}\label{lowlev1d}
Let $A_n=\cir_n(f)$, $P_n(p) = \cir_n(p)K_n^T$, $R_n(r) =
K_n\cir_n(r)$, with $f, p, r$ trigonometric polynomials and $p, r$
satisfying conditions \eqref{NEWipo}. Then
\begin{enumerate}
\item $A_{n/2}=R_n(r)A_nP_n(p)$ coincides with $\cir_{n/2}(\hat{f})$ where
     \begin{equation}\label{eq_fnext}
            \hat{f}(x)=\frac{1}{2^d}\sum_{y\in \Omega(x/2)}r(y)f(y)p(y), \qquad
                x \in [-\pi,\pi)^d.
    \end{equation}
\item If $x^0\in [-\pi,\pi)^d$ is a zero of $f$, then $y^0 = 2x^0 \, \mathrm{mod} \, 2\pi$ is
      a zero of $\hat{f}$. Moreover the order of the zero $y^0$ of $\hat{f}$ is exactly the
      same as the one of the zero $x^0$ of \nolinebreak $f$.
\end{enumerate}
\end{proposition}
\begin{proof}
The essentials of the proof in the case of $r=p$ can be found in
\cite{Serra-Possio}. For $r \neq p$ we can proceed similarly. We
sketch the main steps for the one dimensional case and, at the end,
we extend it to the multidimensional case, mainly for emphasizing
the algebraic interpretation of the frequencies packaging used in
the LFA.

The main relationship is
\begin{equation*}
    K_nF_n=\frac{1}{\sqrt{2}}\left[F_{n/2} \;| \;
    F_{n/2}\right],
\end{equation*}
that implies
\begin{eqnarray*}
  A_{n/2} &=& K_n\cir_n(r)\cir_n(f)\cir_n(p)K_n^T \\
    &=& \frac{1}{2}\left[F_{n/2} \;| \; F_{n/2}\right] \diag_{j=0,\dots,n-1}\left(rfp\left(\frac{2\pi j}{n}\right)\right)
\left[\begin{array}{c}
  F_{n/2}^H \\
  F_{n/2}^H \\
\end{array}\right]\\
&=&\frac{1}{2}F_{n/2}
\diag_{j=0,\dots,n/2-1}\left(rfp\left(\frac{\pi j}{n/2}\right) +
rfp\left(\frac{\pi j}{n/2} + \pi \right)\right)F_{n/2}^H.
\end{eqnarray*}
that is the \eqref{eq_fnext} for $d=1$.

In the multidimensional case $K_nF_n = 2^{-d/2}F_{n/2}G_n$, where
$F_n=\bigotimes_{j=1}^d F_{n_j/2}$ and
$G_n=\bigotimes_{j=1}^d\big([1 \; 1] \otimes I_{n_j/2}\big)$.
Therefore $G_nD_n(rfp)G_n^T = D_{n/2}(2^d\hat{f})$, where
$D_k(h)=\mathrm{diag}_{0\leqslant j \leqslant k-e}(h(2\pi j / n))$
and $\hat{f}$ is defined by (\ref{eq_fnext}).

The claim in item 2 is a consequence of item 1 and of relations
(\ref{NEWipo}).
\end{proof}

Thanks to Proposition \ref{prop:equiv2}, the TGM conditions in
\eqref{NEWipo} give a complete generalization of Proposition
\ref{prop:hemk}, also for non-differential problems since the case
of generating functions vanishing at points different from the
origin is also included. Therefore the analysis based on the zeros
of the generating function can be considered an algebraic
generalization of the LFA also to non-differential problems. For
instance, some discretized integral problems have a generating
function vanishing at $\pi e$, $e=(1,\dots,1)^T \!\in\! \R^d$, or
more generally at some $\pi x$, $x_i \in \{0,1\}$, $i=1,\dots,d$,
$||x||_\infty=1$, with order $2q$. In this case
$\mu_q(x)=2^{-dq}\prod_{j=1}^d(1-\e^{-\imm
x_j})^{\lfloor\frac{q}{2}\rfloor} \e^{\imm x_j
\lfloor\frac{q}{2}\rfloor}$ satisfies the conditions (\ref{NEWipo})
and therefore it defines an optimal TGM. We note that
$R_n=K_{n}T_n(\mu_q)$ is an high-pass filter and then it projects
into the high frequencies. However, thanks to Proposition
\ref{lowlev1d}, the zero at the next level moves to the origin and
at the coarser grids the problem becomes spectrally equivalent to
the discretization of a constant coefficient elliptic PDE.

Finally, we recall the $V$-cycle optimality conditions for Toeplitz
matrices given in \cite{ADS} for $d=1$ and in \cite{ADmultid} for
$d>1$.
\begin{proposition}[\cite{ADS,ADmultid}]\label{V-cycle}
Let $A_{n^{(i)}} = \cir_{n^{(i)}}(f_i)$ be the coefficient matrix at the level $i$,
for $i=0,\dots,l$, with $f_i$ having a unique zero at $x^0_i$. Defining
$A_{n^{(i+1)}} = R_{n^{(i)}}(r_i)A_{n^{(i)}}P_{n^{(i)}}(p_i)$,
$P_{n^{(i)}}(p_i) = \cir_{n^{(i)}}(p_i)K_n^T$ and
$R_{n^{(i)}}(r_i)=\alpha P_{n^{(i)}}(p_i)^H$, i.e. $r_i=\alpha p_i$, $\alpha \in \R \backslash \{0\}$,
where $p_i$ is a trigonometric polynomial non
identically zero and such that for each $x \in [-\pi, \pi)^d$
\begin{subequations}\label{vipo}
    \begin{equation}\label{vipo_a}
        \limsup_{x\to x^0}\left|\frac{p_i(y)}{f_i(x)}\right| = c < +\infty,
        \qquad \forall \, y \in \mirror(x),
    \end{equation}
    \emph{where}
    \begin{equation}\label{vipo_b}
        \sum_{y\in \Omega(x)} p_i(y)^2 > 0.
    \end{equation}
\end{subequations}
Then the $V$-cycle is optimal.
\end{proposition}

We observe that \eqref{vipo_a} defines a stronger condition on the
order of the grid transfer operators with respect to the condition
\eqref{ipo_a}. On the other hand, choosing $c=0$ in \eqref{ipo_a} or
in \eqref{NEWipo_a}, which is equivalent to require that the
\eqref{eq:condF} is satisfied strictly, is usually enough to obtain
an optimal $V$-cycle as numerically shown in
\cite{Serra-Possio,Trot,ADS}. Following the same analysis done for
Proposition \eqref{TGM_Serra}, the Proposition \eqref{V-cycle} could
be generalized to the case of $r\neq p$ and also applied to
non-constant coefficients PDEs. In this last case, condition
\eqref{vipo_a} could be rewrite as
 $\gamma_r+\gamma_p \geq 2m$.

\section{B-spline grid transfer operators}\label{sec:b-spline}
From the discussion at the end of Section \ref{sect:LFA}, the HF is
more important than LF: in fact, for the latter it is enough require
LF$ >0$. A class of grid transfer operators having HF$ =m$ and LF$
>0$ can be defined by
\begin{equation}\label{eq:phi}
    \varphi_m(x)=\prod_{j=1}^d\left(\frac{1+\e^{-\imm x_j}}{2}\right)^m.
\end{equation}
Therefore every grid transfer operator with HF$=m$ has a generating
function of the form $\varphi_m(x)\nu_m(x)$ such that
$\nu_m(x)\neq0$ for all $x \in \mirror(0)$ and $\nu_m(0)=1$. With
$\nu_m \equiv 1$ we obtain a class of projectors with minimal
support for a fixed order $m$. We note that $\sqrt{2}\varphi_m$ is
the symbol of the B-spline of order $m$ in the multiresolution
analysis (MRA) \cite{Choen}. This technique is based on hierarchies
of nested spaces $V_j \subset V_{j+1}$, $j \in \Z$, defined through
a basis generated (for $d=1$) by translations and dilations
$\beta(2^j x -k)$, $k\in \Z$, of a single scaling function $\beta$.
A scaling function satisfies an equation of the type
$\beta(x)=\sqrt{2}\sum_{k \in \Z}h_k\beta(2x-k)$, which expresses
the nestedness of the spaces $V_j={\rm span}\{\beta(2^jx-k),
k\in\Z\}$. Let $\hat{\beta}(v)$ be the Fourier transform of
$\beta(x)$. Then
\begin{equation}\label{eq:four}
    \hat{\beta}(v)=H(v/2)\hat{\beta}(v/2),
\end{equation}
with $H(v)=1/\sqrt{2}\sum_{k\in\Z}h_k\e^{-\imm kv}$. The function
$H(v)$ is called the symbol of $\beta$.
A common application of the MRA is to approximate a high resolution $f \in V_j$
by a coarser function $\tilde{f}\in V_k$ with $k<j$ without losing a lot of information.

In the one dimensional case a simple scaling function is the
Haar-function $\beta(x)=1$ for $x \in [0,1)$ and zero otherwise,
which satisfies the refinement equation $S_1(x) = S_1(2x) +
S_1(2x-1)$. The Haar-function $S_1$ is the simplest B-spline of
order $m=1$, and  $\{S_1(x-k): \, k \in Z\}$ is an orthonormal basis
of $L_2(\R)$. The B-spline of order $m$ can be defined by $S_m = S_1
* S_{m-1}$, where $*$ is the convolution operator. For instance $S_2
= S_1 * S_1$ is such that $S_2(x) =
1/2(S_2(2x)+2S_2(2x-1)+S_2(2x-2))$ and its translated $\tilde
S_2(x):=S_2(x+1)$, $\tilde S_2(x)=1-|x|$ for $x \in [-1,1]$ and zero
otherwise, is known as the hat function. We remark that $S_2 =
2\varphi_2$ and more generally $S_m = 2\varphi_m$.

The functions $S_m$ are not centered, but they can be easy centered
as previously done for $S_2$: instead of $\nu_m \equiv 1$ we take
the shift $\nu_m = \e^{\imm x \lfloor\frac{m}{2}\rfloor}$. In this
way we define a class of centered projectors such that the symbol of
order $m$ is
\begin{equation}\label{eq:classpro}
    \phi_m(x)=\prod_{j=1}^d\left(\frac{1+\e^{-\imm x_j}}{2}\right)^m
    \e^{\imm x \lfloor\frac{m}{2}\rfloor}.
\end{equation}
The $\phi_m$ have HF$ =m$ and LF$ =2$. As previously observed, we
note that a good class of grid transfer operators should have a high
HF, while for the LF the only request is LF $>0$. For $d=1$, the
$\phi_m$, can be obtained using the Tartaglia's triangle as in Table
\ref{tab:triang}. For $d=2$ we take the tensor product of the one
dimensional stencil and so on for $d>2$.
\begin{table}
  \centering
  \caption{The refinement coefficients $h_k \neq 0$, $k\in \Z$ for
    $2^{m-\frac{1}{2}} \phi_m$ in the one dimensional case.}\label{tab:triang}
  \begin{tabular}{|c|c c c c c|}
    \hline
    $m$ & $h_{-2}$ & $h_{-1}$ & $h_0$ & $h_1$ & $h_2$ \\
    \hline
    1 &   & 1 & 1 &   &  \\
    2 &   & 1 & 2 & 1 &  \\
    3 & 1 & 3 & 3 & 1 &   \\
    4 & 1 & 4 & 6 & 4 & 1 \\
    \hline
  \end{tabular}
\end{table}

The class of projectors defined in \eqref{eq:classpro} is a scaled generalization of
$p_{2k}(x)=\sqrt{2}(1+\cos(x))^k$ proposed in \cite{ADS} for even functions
(zeros of order $m=2k$) in the one dimensional case,
indeed $p_{2k} = 2^{1/2-k}\phi_{2k}$. However, scaling factors does not change the
effectiveness of the projector for the Galerkin approach.
We remark that if $m$ is odd than the grid transfer
operators related to $\phi_m$ are not symmetric for vertex centered
discretization, while they are symmetric for cell centered discretization \cite{Yav}.
For instance $[1 \; 3 \; 3 \; 1]$ is the linear interpolation for cell
centered discretization.

We consider for simplicity the one dimensional case, but the
following observations are true also for $d>1$. For vertex centered
discretizations, we are interested in $\varphi_{2k}$. It is easy to
prove that $P_n(\phi_{2})$ is the linear interpolation. Which is the
geometrical meaning of $P_n(\phi_4)$? Which are the relations
between $P_n(\phi_4)$ and $P_n(g_c)$? Answers to these questions
will be given in the next subsection.

\subsection{A quadratic prolongation of order 4} \label{sect:geominterp}
In this subsection we give a geometric interpretation of
$P_n(\phi_4)$ which has HF$ =4$ like the cubic interpolation.

The simplest but useless prolongation is
\begin{prolist}
    \setcounter{ChToep}{0}
    \item $\hspace*{\fill}p(x)\equiv1.\hfill$
\end{prolist}
Without losing in generality we consider $n$ odd.
For $\v{y} \in \C^n$ and $\v{x} \in \C^{\frac{n-1}{2}}$,
$\v{y}=P_n(1)\v{x}=K_n\v{x}$
does not reconstruct constant functions not identically zero.
Particularly, the choice \textbf{(A)} does not provide a good
approximation for the odd components. Therefore, in the standard
MGM, $P_n(p)$ is frequently chosen as the linear interpolation
\begin{prolist}
    \setcounter{ChToep}{1}
    \item  \label{eq:pro2}
           $ \hspace*{\fill}  p(x) = 1+\cos(x). \hfill $
\end{prolist}
\begin{remark}
    The choice (\textbf{B}) compared to the choice (\textbf{A})
    leaves unchanged the even components but reinforces the odd components,
    which are not well approximated by the choice (\textbf{A}), with a linear interpolation.
\end{remark}

When the choice (\textbf{B}) is ineffective, it is usually replaced
with the cubic interpolation. An alternative is given by
$P_n(\phi_4)$ which has HF$=4$ like the cubic interpolation but a
smaller support. This prolongation follows a strategy similar to
that used for deriving choice (\textbf{B}) from (\textbf{A}): it
leaves the linear interpolation for the odd components and
reinforces the even components. From \eqref{eq:classpro}
\[ \phi_4(x)=4^{-d}\prod_{j=1}^d\left(1+\cos(x_j)\right)^2. \]
Therefore, in the one dimensional case, $r(x)=(1+\cos(x))^2/4$ and
\begin{prolist}
    \setcounter{ChToep}{2}
    \item  \label{eq:pro3}
           $ \hspace*{\fill} \D p(x)=(1+\cos(x))^2/2. \hfill $
\end{prolist}
\begin{remark}
    With respect to the choice (\textbf{B}), this choice leaves unchanged the odd components
    but reinforces the even ones with a quadratic approximation:
    \begin{equation}\label{eq:quadint}
        \v{y}_j =
        \left\{
        \begin{array}{l l}
            (\v{x}_k + \v{x}_{k+1})/2, &  j=2k+1, \\
            (\v{x}_{k-1} + 6\v{x}_k + \v{x}_{k+1})/8, & j=2k,
        \end{array}
        \right.
        \qquad k=1,\dots, n,
    \end{equation}
    where we assume $\v{x}_0 = \v{x}_{n+1} = 0$.
\end{remark}

\begin{figure}
    \centering
    \epsfig{file = 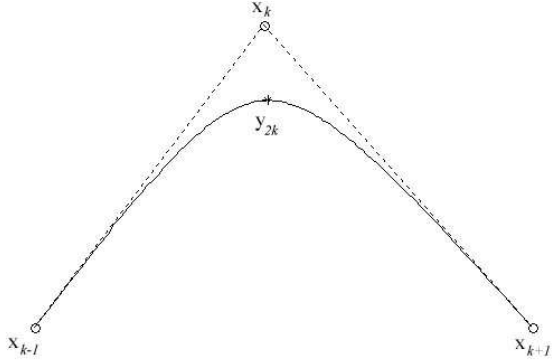, height=5cm}
    \caption{\label{fig:quadpoint} Even components in the finer grid computed with the choice (\textbf{C}).}
\end{figure}

\begin{figure}
    \centering
    \epsfig{file = 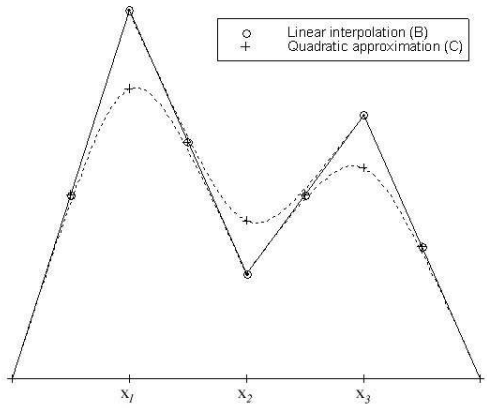, height=6cm}
    \caption{\label{fig:apprfunz} Computation of the points $y_i$ for $i=1,\dots,5$ in
            the finer grid using the linear interpolation (line) and the
            quadratic approximation (dotted).}
\end{figure}

The approximation for the even components of $\v{y}$ is obtained
taking the middle value of a \emph{quadratic rational Bezier curve}
defined from the three points $\{\v{x}_{k-1},
\v{x}_k,\v{x}_{k+1}\}$. The Bernstein polynomial of order $n$ is
defined \nolinebreak as
\[
    B_i^{(n)}(t)=
        \left(  \! \begin{array}{c} n \\  i \\  \end{array} \! \right)
        (1-t)^{i}t^{n-i}, \qquad t \in [0,1], \quad i=1,\dots,n.
\]
A quadratic rational Bezier curve has the expression
\[C(t) = \frac{\sum_{i=0}^{2}\omega_i\v{b}_iB_i^{(2)}(t)}
        {\sum_{i=0}^{2}\omega_iB_i^{(2)}(t)},\]
where $\v{b}_i$ are the control points and $\omega_i$ are the
associated weights for $i=0,1,2$. Let $\v{b}_i = \v{x}_{k+i-1}$ for
$i=0,1,2$, $\omega_1=3/2$ and $\omega_0=\omega_2=1/2$, then
$C(\frac{1}{2})= (\v{x}_{k-1} + 6\v{x}_k + \v{x}_{k+1})/8$ which is
the same of (\ref{eq:quadint}). In Figure \ref{fig:quadpoint} the
previous quadratic approximation is shown for the computation of
$\v{y}_j$ with $j=2k$. Furthermore in Figure \ref{fig:apprfunz} we
compare the values obtained in the finer grid with the choice
(\textbf{B}) (linear interpolation) and with the choice (\textbf{C})
(quadratic approximation).

We note a different philosophy between the cubic interpolation and
this quadratic approximation. The cubic interpolation reconstructs
exactly the polynomial of degree at most three, while the quadratic
approximation does not. However, for the even nodes the cubic
interpolation preserves the exact value as in the coarse grid. This
can be useful for the TGM or when the coarse solution is well
approximated. But when this is not the case, like for the $V$-cycle
in some applications, it could be better to take an approximation
that is an average at the neighboring nodes. The underling idea is
that it is not useful to take a powerful prolongation and a poor
restriction, because in this case we interpolate the solution of a
coarse problem that does not represent well the finer problem.

The generating function of the cubic interpolation is
$g_c(x)=\varphi_4(x)\nu(x)$, with $\nu(x)=2-\cos(x)$. The stencil of
$B_n(g_c)$ is $\frac{1}{32}[-1 \quad 0 \quad 9 \quad 16 \quad 9
\quad 0 \quad -1]$ and the one of $B_n(\phi_{4})$ is $\frac{1}{16}[1
\quad 4 \quad 6 \quad 4 \quad 1]$. In addition to the different
philosophy previously emphasized, we observe the following mainly
differences.
\begin{remark}\label{rem1}
  $B_n(g_c)$ has HF=LF=4, while $B_n(\phi_{4})$ has HF=4 and LF=2.
\end{remark}
\begin{remark}\label{rem2}
    From a computational point of view,
    the two stencils have the same number of nonzero elements.
    Hence they have the same computational cost for the projection of a vector
    between coarse to fine or fine to coarse grids.
    The main difference is that the stencil of the cubic interpolation has a lager support.
    This implies that, mainly for $d>1$, $B_n(\phi_{4})$ defines a MGM
    that is computationally more efficient at the coarser grids.
    Indeed, using the Galerkin approach, the coarse matrix has a lower
    bandwidth. Moreover, $B_n(\phi_{4})$ requires less boundary points,
    property that is very useful for a parallel implementation. In
    this way less communications are required among the nodes and we can employ
    a colored Gauss-Seidel with a smaller number of colors, increasing the parallelism degree.
\end{remark}

From the previous remarks, as it will be confirmed by the numerical
experimentation in Section \ref{sect:exp}, the cubic interpolation
usually converges within less iterations with respect to the choice
(\textbf{C}) (see Remark \ref{rem1}), but both have the same
asymptotic behavior since they have the same HF. Therefore, thanks
to Remark \ref{rem2}, $B_n(\phi_{4})$ could be a good alternative to
$B_n(g_c)$, mainly for parallel implementations or $V$-cycle MGMs.

\section{A comparison between MGMs and wavelets methods} \label{sec:wave}
In this subsection we recall some approximation properties of the $\phi_m$
defined in \eqref{eq:classpro}, or equivalently of the $\varphi_m$ in \eqref{eq:phi}.
Moreover we discuss some relations between  wavelets and multigrid methods.

Unfortunately $\{S_m(x-k): \, k\in\Z, m > 1\}$ is not an orthogonal system.
However, the $S_m$ satisfy the quasi-interpolant Strang-Fix condition \cite{Choen}:
\begin{equation}\label{Str-Fx}
\left(\frac{\partial}{\partial v}\right)^s \hat S_m(2k\pi) = 0,
\qquad k\in \Z\backslash\{0\}, \quad s=0,\dots, m-1.
\end{equation}
It follows that the polynomials of degree at most $m-1$ are
contained in the space $V_0= {\rm span}\{S_m(x-k): \, k\in\Z\}$. The
same property is also usually expressed in terms of vanishing
moments. Let $\psi_m$ be the wavelet associated to $S_m$. The first
$m-1$ moments of $\psi_m$ vanish, i.e.
\begin{equation}\label{eq:moments}
    \int_{-\infty}^{+\infty}{x^s\psi_m(x)}=0,\qquad s=1,\dots,m-1.
\end{equation}
Starting from the orthogonality condition
\begin{equation}\label{eq:ortho}
    |H(v)|^2+|H(v+\pi)|^2=1,
\end{equation}
where $H(v)$ is the symbol defined in \eqref{eq:four}, and imposing
the vanishing of the moments, Daubechies defined orthogonal wavelets
\cite{Dau}. The moment condition \eqref{eq:moments} is equivalent to
require $H(x)=\varphi_m(x)\nu(x)$, such that $\nu(x)$ is a
trigonometric polynomial with $\nu(0)=1$. This means that orthogonal
wavelets are obtained imposing orthogonality to $S_m$. For instance,
for $m=2$ we obtain the Daubechies wavelet of order 2 with scaling
coefficients $h=\frac{1}{4\sqrt{2}}[1+\sqrt{3}, \quad 3+\sqrt{3},
\quad 3-\sqrt{3}, \quad 1-\sqrt{3}]$ and symbol
$\varphi_2(x)(1+\sqrt{3} + (1-\sqrt{3}))\e^{-\imm x})/2$.

Considering transfer grid operators for multigrid methods, a
question now seems to be natural: ``Is the orthogonality condition
\eqref{eq:ortho} necessary or the moment condition
\eqref{eq:moments} is sufficient?'' The answer is that the
orthogonality condition \eqref{eq:ortho} is not necessary. Some
motivations were firstly given in \cite{Briggs}. Indeed, in the MRA
we would have $V_0=V_1 \oplus W_1$, while for multigrid methods we
would have $\Omega={\rm Range}\{P\} \oplus {\rm NullSpace}\{R\}$
where $\Omega^h$ is the fine grid \cite{Tutorial}. The error $e$ in
$\Omega^h$ can be decomposed as $e=s+t$, with $s \in {\rm
Range}\{P\}$ and $t \in {\rm NullSpace}\{P^T\}$ (we consider the
Galerkin condition $R=P^T$). We note that $CGCs=0$, while $CGCt=t$.
Therefore, if after the pre-smoothing step $t=0$, then the TGM
converges in one iteration. Clearly this is a too strong condition.
However ${\rm Range}\{P\}$ is spanned by smooth functions while
${\rm NullSpace}\{P^T\}$ is spanned by oscillating functions, thus
if we have a ``good'' smoother $t \approx 0$. Consequently, to
define a good grid transfer operator we are interested only to the
scaling function and we do not use the wavelet since the smoother
has already reduced the error in the high frequencies. Eventually,
the multigrid methods is an iterative method, hence it is not
necessary to have convergence in one iteration but the solution can
be substantially improved iterating. These considerations show that
for multigrid methods the orthogonality condition \eqref{eq:ortho}
is not necessary. The only necessary condition concerns moments in
\eqref{eq:moments}, that is equivalent to the HF and to the
factorization $\varphi_m\nu_m$.

Concluding, for Toeplitz linear system relations
between wavelets and multigrid methods was already investigated in \cite{TGM-wave}.
In such paper the authors proved the TGM optimality using Cohen, Daubechies and Feauveau (CDF) 9/7
biorthogonal wavelets for generating functions having a zero of order 4 in the origin.
Thanks to the previous comments this is expected since CDF 9/7 have both moments of order 4.
However the interesting fact is that the proof in \cite{TGM-wave} is done directly in the
Toeplitz class and not in the algebra case. Moreover, computing a class of compactly supported
biorthogonal wavelet systems GCDF with specified vanishing moments for scaling functions,
the authors show numerically that for moments of order $k$ the TGM is optimal for
generating functions having a zero of order at most $2k$.
This is exactly the condition \eqref{ipo_a} since the symbols have the form
$H(x)=\varphi(x)\nu(x)$.

\section{Numerical Experiments}\label{sect:exp}
For the sake of simplicity, we consider some tests in the 1D case.
The results in the multidimensional case are very similar. The
numerical experiments are done using Matlab 7.0. We fix the
following parameters. The MGMs are stopped when the relative $l_2$
norm of residual is lower than $10^{-9}$. The coarse matrices are
defined using the Galerkin approach. The coarser problem has
dimensions $7\times 7$. These means that
$A_{n^{(i)}}=T_{n^{(i)}}(f_i) \;+$ low rank, for $i=1,\dots,l$,
where the $f_i$ are defined according to~\eqref{eq_fnext}. The right
hand side is obtained from the exact solution $x_j = j/n$,
$j=1,\dots,n$. The pre-smoother and the post-smoother are one step
of relaxed Richardson with relaxation parameter equal to
$1.5/\norm[\infty]{f_i}$ and $1/\norm[\infty]{f_i}$ according to
\cite{ADmultid}.

In Section \ref{sec:gentoep}, Proposition \ref{prop:equiv2} gives a
validation of the TGM conditions \eqref{NEWipo} for PDEs. Here, in
the first test problem, we give a numerical validation of the
conditions \eqref{NEWipo} in the case of discretized integral
problems. We show that sometimes the choice $r\neq p$  can be very
useful also for the MGM for Toeplitz matrices. In the second test
problem, we show that the $V$-cycle optimality result for Toeplitz
matrices in Proposition \ref{V-cycle} can be applied also to
non-constant coefficients PDEs.

\subsection{An integral problem}
We consider the discretization by the rectangle quadrature formula
of the following Fredholm operator of first kind:
\[
g(x)\, =\,\int_{{\R}}{t(x-\theta) \, f(\theta)\, \diff \theta},
\qquad x,\, \theta \in \R,
\]
where $f$ is the input object, $t$ is the integral kernel of the
operator, also called point spread function (PSF) and $g$ is the
observed object. In the discrete case, when zero Dirichlet BCs are
used and the PSF is shift invariant, the above approximation gives
rise to the system $A_n{\bf f}={\bf g}$, where $A_n=\toe_n(z)$ with
$z(\pi)=0$ and positive elsewhere. According to the analysis in
Section \ref{sec:gentoep}, the generating function of a grid
transfer operator of order $s$ will be
 $\mu_s$ at the finer level and $\phi_s$ at the lower level.
We denote by $\delta_r$ and $\delta_p$ the order of the zero of the
restriction and of the prolongation, respectively.

First of all, we give a numerical evidence of the relevance of
conditions \eqref{NEWipo}. We consider $z(x)=(2+2\cos(x))^3$, which
has a zero in $\pi$ of order 6. From Proposition \ref{TGM_Serra}
$\delta_r=\delta_p=2$ is not enough for an optimal TGM, while it is
necessary to set $\delta_r=\delta_p=4$. However, if we allow
$\delta_r \neq \delta_p$ then, from conditions \eqref{NEWipo}, the
optimality of the TGM is guaranteed with $\delta_r=2$ and
$\delta_p=4$. This is confirmed by the numerical results in Table
\ref{tab:ex1tgm}. In such table we report the number of iterations
required by the TGM for converging for different orders of the grid
transfer operators, when increasing the problem size.
\begin{table}
   \centering
    \caption{\label{tab:ex1tgm} TGM iteration numbers varying the order
        of the grid transfer operators and the problem size $n$
        for the integral problem $z(x)=(2+2\cos(x))^3$. }
   \begin{tabular}{| c | c  c  c |}\hline
         $n$ & $\delta_r=2$ & $\delta_r=2$ & $\delta_r=4$\\
             & $\delta_p=2$ & $\delta_p=4$ & $\delta_p=4$ \\   \hline
        $15$  &  219  & 65 & 51  \\
        $31$  &  607  & 72 & 52  \\
        $63$  &  1501  & 76 & 51 \\
        $127$ &  $>$ 2000  & 77 & 50 \\
        $255$ &  $>$ 2000  & 78 & 49 \\
        \hline
    \end{tabular}
\end{table}

In real applications usually more than two-grids are used. The TGM
optimality conditions are not enough for the $V$-cycle, but they
give good estimations for the $W$-cycle \cite{Trot}. Defining
$z(x)=(2+2\cos(x))^2$, which has a zero in $\pi$ of order 4, also
the choice $\delta_r=\delta_p=2$ gives an optimal TGM. However, in
Table \ref{tab:ex1w} we see that for the $W$-cycle there is a large
reduction in the iteration number for $\delta_r=2$ and $\delta_p=4$,
with respect to $\delta_r=\delta_p=2$. Furthermore using
$\delta_r=\delta_p=4$ implies a negligible reduction of the
iteration number, with respect to the choice $\delta_r=2$ and
$\delta_p=4$.

From a computational point of view, we should investigate the
structure of the coefficient matrices at each level for the previous
choices of the grid transfer operators. Indeed, at each level the
main computational cost is related to the matrix vector product with
the matrices $A_{n^{(i)}}$, for $i=0,\dots, l$. In the last example
with $z(x)=(2+2\cos(x))^2$, we use the following grid transfer
operators of order $2s$: $r_0=p_0=(2-2\cos(x))^s$  and
$r_i=p_i=(2+2\cos(x))^s$ for $i=1,\dots,l-1$. For
$\delta_r=\delta_p=2$ we have $A_{n^{(i)}}=T_n^{(i)}(\tilde z)$, for
$i=1,\dots, l$, where $\tilde z(x)=(2-2\cos(x))^2$. For $\delta_r=2$
and $\delta_p=4$ we have $A_{n^{(i)}}=2^iT_n^{(i)}(\tilde
z)+c_ie_1e_1^T + c_ie_ne_n^T$, where $e_j$ is the $j$th vector of
the canonical basis and $c_i =
[T_n^{(i)}(r_i)A_{n^{(i)}}T_n^{(i)}(p_i)]_{2,2}-2^i\cdot 6$.
Therefore, the matrix vector product has about the same
computational cost for both choices. On the other hand, for
$\delta_r=\delta_p=4$ the matrix $A_{n^{(i)}}$ is Toeplitz plus a 4
rank correction and moreover the bandwidth of the Toeplitz part is
not longer 5, but it becomes 7. Obviously, this fact increases the
complexity and the computational cost of the matrix vector product.
Moreover, for using the Gauss-Seidel smoother with a coloring
strategy we need 4 colors instead of 3 colors as for $T_n(z)$,
losing a degree of parallelism. The previous considerations are
enhanced in the multidimensional case. Indeed in the two dimensional
case the bandwidth of each block moves from 5 to 7, and also the
block bandwidth moves from 5 to 7.

For preserving the Toeplitz structure at each level, a different
cutting matrix proposed in \cite{ADS} can be used. The main idea is
in the changing of the coarse problems size in order to neglect in
some way the boundary effects that give the low rank corrections.
However, in this way we lose some information and the iteration
number slightly increases even if the general (optimal) behavior is
preserved in the numerical experimentations.

\begin{table}
   \centering
   \begin{footnotesize}
    \caption{\label{tab:ex1w} $W$-cycle iteration numbers varying the order
        of the grid transfer operators and the problem size $n$
        for the integral problem $z(x)=(2+2\cos(x))^2$.}
   \begin{tabular}{| c | c  c  c |}  \hline
         $n$ & $\delta_r=2$ & $\delta_r=2$ & $\delta_r=4$ \\
             & $\delta_p=2$ & $\delta_p=4$ & $\delta_p=4$ \\   \hline
        $31$  & 25 & 23 & 22 \\
        $63$  & 32 & 23 & 21 \\
        $127$ & 35 & 23 & 21 \\
        $255$ & 37 & 23 & 20 \\
        $511$ & 37 & 23 & 20 \\
        \hline
    \end{tabular}
   \end{footnotesize}
\end{table}

\subsection{A differential equation with nonconstant coefficients}
We consider the following equation
\begin{equation}\label{eq:der4}
    \left\{
    \begin{array}{l}\displaystyle
        \frac{\diff^2}{\diff x^2} \left(a(x) \frac{\diff^{2}}{\diff x^2}u(x)\right) = g(x), \qquad x \in (0,1),\\
        u(0)=u(1)=0
    \end{array}\right.
\end{equation}
with nonconstant $a(x)$ and order $m=4$.

In this subsection we consider the $V$-cycle, that is cheaper than
the $W$-cycle in vector and parallel implementations. Therefore we
need a more powerful smoother and we replace the weighted Richardson
with Gauss-Seidel. The generating function of the grid transfer
operator is the same for each coarse problem. We test several
combinations of $\phi_k$, $k=2,4$ and the cubic interpolation $g_c$.

\begin{table}
   \centering
   \begin{footnotesize}
    \caption{\label{tab:der4} $V$-cycle iteration numbers varying problem size $n$
         and $a(x)$ for the differential problem \eqref{eq:der4} ($g_c$ = cubic interpolation).}
   \begin{tabular}{| c | c  c  c  c c | c  c  c  c c |}  \hline
         restriction & $\phi_2$ & $\phi_2$ & $\phi_2$ & $\phi_4$ & $\phi_4$
                    & $\phi_2$ & $\phi_2$ & $\phi_2$ & $\phi_4$ & $\phi_4$\\
         prolongation & $\phi_2$ & $\phi_4$ & $g_c$ & $\phi_4$ & $g_c$
                    & $\phi_2$ & $\phi_4$ & $g_c$ & $\phi_4$ & $g_c$ \\   \hline
         $n$  & \multicolumn{5}{c |}{$a(x) = \e^x$}
             & \multicolumn{5}{c |}{$a(x) = (x-0.5)^2$} \\ \hline
        15 & 14 & 9 & 9 & 7 & 7 &15 & 10   & 10  & 9  & 9\\
        31 & 32 & 11 & 13& 10 & 9 & 33 & 13  & 17  & 10 & 11\\
        63 & 60 & 17 & 15& 14 & 9 & 61 & 17  & 24  & 13 & 11\\
        127 & 98 & 27 & 20& 18 & 12 & 101 & 26 & 27  & 17 & 13\\
        255 & 151& 38& 27& 22& 16 & 155 & 35 & 29  & 20 & 16\\
        511 & 215 &48 & 34 &26& 20 & 221 & 44 & 36  & 24 & 19\\
        1023 & 276& 57& 44&29& 22 & 284 & 53 & 46 & 27 & 22\\
        \hline
    \end{tabular}
   \end{footnotesize}
\end{table}

We note that in practical implementations it is usually required
that \eqref{eq:condF} is satisfied strictly. This condition in terms
of generating functions is equivalent to require that $c=0$ in
\eqref{NEWipo_a}: the latter is a generalization of a similar
condition already numerically observed in \cite{ADS} for $r=p$ and
when considering the $V$-cycle. Indeed, from Table \ref{tab:der4},
we see that the choice $(\phi_2, \phi_2)$ for the couple
(restriction, prolongation) is not effective. We recall that
$\phi_2$ is the linear interpolation. For obtaining an effective
$V$-cycle it is enough to increase only the order of the
prolongation (i.e., $\gamma_r = 2$ and $\gamma_p = 4$). However,
more stable results can be obtained for $\gamma_r = 4$ and $\gamma_p
= 4$ according to Proposition \ref{V-cycle}.

Eventually, we compare the two prolongations generated by $\phi_4$
and $g_c$. In terms of number of iterations the function $g_c$ has
to be preferred with respect to $\phi_4$ according to Remark
\ref{rem1}. On the other hand, $\phi_4$ has a 5-point stencil while
$g_c$ has a 7-point stencil. This means that the choice ($\phi_2,
g_c$) leads to coarse matrices having a bandwidth equal to 7 like
the choice ($\phi_4, \phi_4$), while ($\phi_4, g_c$) leads to a
bandwidth equal to 9. Therefore, the right comparison should be
between ($\phi_2, g_c$) and ($\phi_4, \phi_4$). From Table
\ref{tab:der4}, it is evident that the second one has to be
preferred according to Remark \ref{rem2}. Moreover, considering the
discretization near the boundary, $\phi_k$ requires less boundary
points than $g_c$.

\section{Conclusions}\label{sect:concl}
Considering elliptic PDEs with constant coefficients, we have shown
the equivalence between the LFA and the analysis based on the zeros
of the generating functions of Toeplitz matrices. This equivalence
has two implications. The first one is that the techniques used for
Toeplitz and Circulant matrices allow to extend the LFA also to
non-differential problems (e.g., integral problems). The second one
is that it suggests to choose the restriction different from the
prolongation also for MGMs for Toeplitz linear systems. The
generalization of the MGM for Toeplitz matrices proposed in Section
\ref{sec:gentoep} replaces the Galerkin conditions with the
following:
\begin{enumerate}
  \item $A_{n/2} =  R_{n}(r) A_{n} P_{n}(p)$
  \item $r$ not necessary equal to $p$, but such that $p \geq 0$,
  $r \geq 0$ and both even or odd (such that $A_{n/2}$ is again positive definite).
\end{enumerate}

We have given a class of grid transfer operator with minimal support
for a fixed HF and the geometrical interpretation of the operator
with HF $=4$. Such class is related to the symbol of B-spline,
giving the possibility to discuss some useful and suggestive
relations between wavelets and multigrid methods.


\end{document}